\documentclass[a4paper]{article}
\usepackage{amsmath,amssymb,mathrsfs,graphicx,psfrag}
\DeclareMathOperator{\re}{Re}
\DeclareMathOperator{\im}{Im}
\DeclareMathOperator{\fp}{f.p.}
\newtheorem{thm}{Theorem}
\allowdisplaybreaks[1]
\begin{document}
\title{A numerical method for computing Hadamard finite-part integrals with a non-integral power singularity
at an endpoint}
\author{Hidenori Ogata%
\footnote{%
Department of Computer and Network Engineering, 
Graduate School of Informatics and Engineering, 
The University of Electro-Communications, 
1-5-1 Chofugaoka, Chofu, Tokyo 182-8585, Japan, 
(e-mail) {\tt ogata@im.uec.ac.jp}
}}
\maketitle
\begin{abstract}
 In this paper, we propose a numerical method of computing a Hadamard finite-part integral with a non-integral power 
 singularity at an endpoint, that is, a finite part of a divergent integral as a limiting procedure. 
 In the proposed method, we express the desired finite-part integral using a complex loop integral, and 
 obtain the finite-part integral by evaluating the complex integral by the trapezoidal formula. 
 Theoretical error estimate and some numerical examples show the effectiveness of the proposed method. 
\end{abstract}
\section{Introduction}
\label{sec:introduction}
The integral
\begin{equation*}
 \int_0^1 x^{\alpha-2}f(x)\mathrm{d}x \quad ( \: 0 < \alpha < 1 \: ),
\end{equation*}
where $f(x)$ is an analytic function on the closed interval $[0,1]$, is divergent. 
However, we can assign a finite value to this divergent integral as follows. 
For $0 < \epsilon < 1$, using integration by part, we have
\begin{align*}
 \int_{\epsilon}^1 x^{\alpha-2}f(x)\mathrm{d}x 
 = \: & 
 \frac{1}{\alpha-1}\int_{\epsilon}^{1} (x^{\alpha-1})^{\prime}f(x)\mathrm{d}x
 \\ 
 = \: & 
 \frac{1}{\alpha-1}
 \bigg[x^{\alpha-1}f(x)\bigg]_{\epsilon}^{1}
 - 
 \int_{\epsilon}^{1}x^{\alpha-1}f^{\prime}(x)\mathrm{d}x
 \\ 
 = \: & 
 \frac{\epsilon^{\alpha-1}f(\epsilon)-f(1)}{1-\alpha} 
 - 
 \int_{\epsilon}^{1}x^{\alpha-1}f^{\prime}(x)\mathrm{d}x
 \\ 
 = \: & 
 \frac{\epsilon^{\alpha-1}f(0)}{1-\alpha} + 
 \mbox{(terms finite as $\epsilon\downarrow 0$)}, 
\end{align*}
and the limit
\begin{equation*}
 \lim_{\epsilon\downarrow 0}
  \left\{
   \int_{\epsilon}^{1}x^{\alpha-2}f(x)\mathrm{d}x 
   - 
   \frac{\epsilon^{\alpha-1}f(0)}{1-\alpha}
  \right\}
\end{equation*}
is finite. 
We call this limit an Hadamard finite-part (f.p.) integral and denote it by
\begin{equation*}
 \fp\int_0^{1} x^{\alpha-2}f(x)\mathrm{d}x.
\end{equation*}
Similarly, we can define a f.p. integral
\begin{equation}
 \label{eq:fp-integral0}
  \fp\int_0^1 x^{\alpha-1-n}f(x)\mathrm{d}x
\end{equation}
for $n=1, 2, \ldots$, $0 < \alpha < 1$ and an analytic function $f(x)$ on $[0,1]$ \cite{EstradaKanwal1989}. 

In this paper, we propose a numerical method of computing f.p. integrals (\ref{eq:fp-integral0}). 
In the proposed method, we express the f.p. integral using a complex loop integral, and we obtain 
the desired f.p. integral by evaluating the complex integral by the trapezoidal formula with equal mesh. 
Theoretical error estimate and numerical examples will show that the approximation by the proposed method 
converges exponentially as the number of sampling points increases.

Previous works related to this paper are as follows. 
The author and Hirayama proposed a numerical integration method for ordinary integrals related to hyperfunction theory 
\cite{OgataHirayama2018}, 
where a desired integral is expressed using a complex loop integral, and it is obtained by evaluating the complex 
integral by the trapezoidal formula with equal mesh.  
The author proposed a numerical method of computing f.p. integrals 
with an integral order singularity \cite{Ogata2019b}.
For Cauchy principal value integrals and Hadamard finite-part integrals with a singularity 
in the interior of the integral interval 
\begin{equation}
 \label{eq:fp-integral02}
  \fp\int_0^1 \frac{f(x)}{(x-\lambda)^n}\mathrm{d}x \quad ( \: 0 < \lambda < 1, \: n = 1, 2, \ldots \: ), 
\end{equation}
many numerical methods were proposed. 
Elliot and Paget proposed a Gauss-type numerical integration formulas for f.p. integrals 
(\ref{eq:fp-integral02}) \cite{ElliotPaget1979,Paget1981}. 
Bialecki proposed Sinc numerical integration formulas for f.p. integrals \cite{Bialecki1990a,Bialecki1990b}, 
where the trapezoidal formula with the variable transform technique are used as in the DE formula for 
ordinary integrals \cite{TakahasiMori1974}. 
The author et al. improved them and proposed a DE-type numerical integration formulas for f.p. integrals 
(\ref{eq:fp-integral02}) \cite{OgataSugiharaMori2000}. 

The remainder of this paper is structured as follows. 
In Section \ref{sec:fp-integral}, we define the f.p. integrals and propose a numerical method of computing them. 
Then, we give a theorem on error estimate of the proposed method. 
In Section \ref{sec:example}, we show some numerical examples which show the 
effectiveness of the proposed method. 
In Section \ref{sec:summary}, we give a summary of this paper. 

\section{Hadamard finite-part integral}
\label{sec:fp-integral}
The Hadamard finite-part integral is defined by 
\begin{multline}
 \label{eq:fp-integral}
 \fp\int_0^1 x^{\alpha-1-n}f(x)\mathrm{d}x = 
 \lim_{\epsilon\downarrow 0}
 \left\{
 \int_{\epsilon}^1 x^{\alpha-1-n}f(x)\mathrm{d}x
 - 
 \sum_{k=0}^{n-1}\frac{\epsilon^{\alpha-n+k}}{k!(n-k-\alpha)}f^{(k)}(0)
 \right\}
 \\ 
 ( \: n = 1, 2, \ldots; \: 0 < \alpha < 1 \: ).
\end{multline}
We can see that it is well-defined using integration by part. 
In fact, repeating integration by part, we have
\begin{align*}
 & 
 \int_{\epsilon}^1 x^{\alpha-n-1}f(x)\mathrm{d}x
 \\ 
 = \: & 
 \frac{\epsilon^{\alpha-n}}{n-\alpha}f(\epsilon) + 
 \frac{\epsilon^{\alpha-n+1}}{(n-\alpha)(n-1-\alpha)}f^{\prime}(\epsilon) + 
 \frac{\epsilon^{\alpha-n+2}}{(n-\alpha)(n-1-\alpha)(n-2-\alpha)}f^{\prime\prime}(\epsilon)
 \\ 
 & 
 + \cdots + 
 \frac{\epsilon^{\alpha-1}}{(n-\alpha)(n-1-\alpha)(n-2-\alpha)\cdots(1-\alpha)}f^{(n-1)}(0)
 \\
 & 
 + 
 \frac{1}{(n-\alpha)(n-1-\alpha)(n-2-\alpha)\cdots(1-\alpha)}
 \int_{\epsilon}^{1}x^{\alpha-1}f^{(n)}(x)\mathrm{d}x
 \\
 & 
 + 
 \mbox{(terms finite as $n\downarrow 0$, which will be denoted by ^^ ^^ $\cdots$'' below)}
 \\ 
 = \: & 
 \frac{\epsilon^{\alpha-n}}{n-\alpha}
 \sum_{k=0}^{n-1}\frac{\epsilon^k}{k!}f^{(k)}(0)
 + 
 \frac{\epsilon^{\alpha-n+1}}{(n-\alpha)(n-1-\alpha)}
 \sum_{k=0}^{n-2}\frac{\epsilon^k}{k!}f^{(k)}(0)
 \\ 
 & 
 + \frac{\epsilon^{\alpha-n+2}}{(n-\alpha)(n-1-\alpha)(n-2-\alpha)}
 \sum_{k=0}^{n-3}\frac{\epsilon^k}{k!}f^{(k)}(0)
 \\ 
 & 
 + \cdots + 
 \frac{\epsilon^{\alpha-1}}{(n-\alpha)(n-1-\alpha)(n-2-\alpha)\cdots(1-\alpha)}f^{(n-1)}(0)
 + \cdots 
 \\ 
 = \: & 
 \frac{\epsilon^{\alpha-n}}{n-\alpha}f(0) + 
 \frac{\epsilon^{\alpha-n+1}}{n-\alpha}\left(1 + \frac{1}{n-1-\alpha}\right)f^{\prime}(0)
 \\ 
 & 
 + 
 \frac{\epsilon^{\alpha-n+2}}{n-\alpha}
 \left\{
 \frac{1}{2!} + \frac{1}{n-1-\alpha} + \frac{1}{(n-1-\alpha)(n-2-\alpha)}
 \right\}
 f^{\prime\prime}(0)
 \\ 
 & 
 + 
 \frac{\epsilon^{\alpha-n+3}}{n-\alpha}
 \left\{
 \frac{1}{3!} + \frac{1}{2!(n-1-\alpha)} + \frac{1}{(n-1-\alpha)(n-2-\alpha)}  
 \right.
 \\
 & 
 \hspace{18mm}
 +
 \left.
 \frac{1}{(n-1-\alpha)(n-2-\alpha)(n-3-\alpha)}
 \right\}
 f^{\prime\prime\prime}(0)
 \\ 
 & 
 + \cdots  
 \\ 
 & + 
 \frac{\epsilon^{\alpha-1}}{n-\alpha}
 \left\{
 \frac{1}{(n-1)!(n-1-\alpha)} + \frac{1}{(n-2)!(n-1\alpha)(n-2-\alpha)}  
 \right.
 \\
 & 
 \hspace{18mm}
 +
 \frac{1}{(n-3)!(n-1-\alpha)(n-2-\alpha)(n-3-\alpha)} 
 \\
 & 
 \hspace{18mm}
 \left.
 + \cdots + 
 \frac{1}{(n-1-\alpha)(n-2-\alpha)\cdots(1-\alpha)}
 \right\}
 f^{(n-1)}(0) 
 \\ 
 & + \cdots 
 \\ 
 = \: & 
 \frac{\epsilon^{\alpha-n}}{n-\alpha}f(0) + 
 \frac{\epsilon^{\alpha-n+1}}{n-1-\alpha}f^{\prime}(0)
 + 
 \frac{\epsilon^{\alpha-n+2}}{2!(n-2-\alpha)}f^{\prime\prime}(0)
 + \cdots + 
 \frac{\epsilon^{\alpha-1}}{(n-1)!(1-\alpha)}f^{(n-1)}(0)
 \\ 
 & 
 + \cdots.
\end{align*}

\begin{figure}[htbp]
 \begin{center}
  \psfrag{0}{$\mathrm{O}$}
  \psfrag{1}{$1$}
  \psfrag{p}{$\Gamma_{\epsilon}^{(+)}$}
  \psfrag{m}{$\Gamma_{\epsilon}^{(-)}$}
  \psfrag{a}{$C_{\epsilon}^{(0)}$}
  \psfrag{b}{$C_{\epsilon}^{(1)}$}
  \psfrag{e}{$\epsilon$}
  \psfrag{f}{$(1+\epsilon)^{-1}$}
  \includegraphics[width=0.6\textwidth]{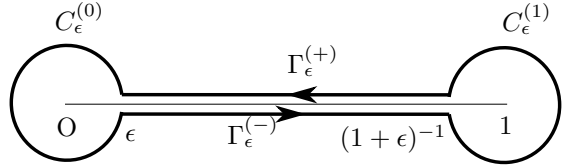}
 \end{center}
 \caption{The integral paths.}
 \label{fig:proof-integral-path}
\end{figure}

If the integrand $f(x)$ is an analytic function on the closed interval $[0,1]$, 
the f.p. integral (\ref{eq:fp-integral}) is expressed using a complex loop integral as in the following 
theorem.
\begin{thm}
 \label{thm:complex-integral}
 We suppose that $f(z)$ is an analytic function in a complex domain $D$ 
 containing the closed interval $[0,1]$ in its interior. 
 Then, the f.p. integral (\ref{eq:fp-integral}) is expressed as 
 \begin{gather}
  \label{eq:complex-integral}
  \fp\int_0^1 x^{\alpha-1-n}f(x)\mathrm{d}x 
  = 
  \frac{1}{2\pi\mathrm{i}}\oint_C z^{-n}f(z)\Psi_{\alpha}(z)\mathrm{d}z 
  + 
  \sum_{k=0}^{n-1}\frac{f^{(k)}(0)}{k!(\alpha-n+k)}, 
  \intertext{where}
  \label{eq:Psi}
  \Psi_{\alpha}(z) = 
  \alpha^{-1}z^{-1}F(\alpha, 1; \alpha+1; z^{-1}), 
 \end{gather}
 and $C$ is a closed complex integral path contained in $D\setminus[0,1]$ 
 and encircling the interval $[0,1]$ in the positive sense. 
\end{thm}
\paragraph{Proof of Theorem \ref{thm:complex-integral}}
From Cauchy's integral theorem, the complex integral of the first term 
on the right-hand side of (\ref{eq:complex-integral}) is modified into 
\begin{equation*}
 \frac{1}{2\pi\mathrm{i}}\oint_C z^{-n}f(z)\Psi_{\alpha}(z)\mathrm{d}z
  = 
  \left(
   \int_{C_{\epsilon}^{(0)}} + \int_{\Gamma_{\epsilon}^{(-)}}
   + \int_{C_{\epsilon}^{(1)}} + \int_{\Gamma_{\epsilon}^{(+)}}
  \right)
  z^{-n}f(z)\Psi_{\alpha}(z)\mathrm{d}z,
\end{equation*}
where $C_{\epsilon}^{(0)}$, $C_{\epsilon}^{(1)}$, $\Gamma_{\epsilon}^{(+)}$ and 
$\Gamma_{\epsilon}^{(-)}$ are complex paths respectively defined by 
\begin{align*}
 C_{\epsilon}^{(0)} = \: & 
 \{ \: \epsilon\mathrm{e}^{\mathrm{i}\theta} \: | \: 0 \leqq \theta \leqq 2\pi \: \}, 
 \\ 
 C_{\epsilon}^{(1)} = \: & 
 \{ \: (1-\epsilon\mathrm{e}^{\mathrm{i}\theta})^{-1} \: | \: 0 \leqq \theta \leqq 2\pi \: \}, 
 \\ 
 \Gamma_{\epsilon}^{(+)} = \: & 
 \{ \: x\in\mathbb{R}+\mathrm{i}0 \: | \: 1-\epsilon \geqq x \geqq \epsilon \: \}
 \\
 \Gamma_{\epsilon}^{(-)} = \: & 
 \{ \: x\in\mathbb{R}-\mathrm{i}0 \: | \: \epsilon \leqq x \leqq 1 - \epsilon \: \}
\end{align*}
with $0 < \epsilon \ll 1$ (see Figure \ref{fig:proof-integral-path}). 
From the formula 15.3.7 in \cite{AbramowitzStegun1965}, we have
\begin{equation*}
 \Psi_{\alpha}(z) = 
 \frac{-\pi}{\sin\pi\alpha}(-z)^{\alpha-1} - \frac{1}{\alpha-1}F(1-\alpha,1;2-\alpha;z)
 \quad ( \: |\arg(-z)|<\pi \: ). 
\end{equation*}
Then, as to the integrals on $\Gamma_{\epsilon}^{(\pm)}$, we have 
\begin{align*}
 & 
 \frac{1}{2\pi\mathrm{i}}\left( \int_{\Gamma_{\epsilon}^{(+)}} + \int_{\Gamma_{\epsilon}^{(1)}}\right)
 z^{-n}f(z)\Psi_{\alpha}(z)\mathrm{d}z
 \\
 = \: & 
 \frac{-1}{2\mathrm{i}\sin\pi\alpha}
 \left\{
 - \int_{\epsilon}^{1-\epsilon}x^{-n}f(x)(-(x+\mathrm{i}0))^{\alpha-1}\mathrm{d}x
 + \int_{\epsilon}^{1-\epsilon}x^{-n}f(x)(-(x-\mathrm{i}0))^{\alpha-1}\mathrm{d}x
 \right\}
 \\
 = \: & 
 \frac{-1}{2\mathrm{i}\sin\pi\alpha}
 (-\mathrm{e}^{-\mathrm{i}\pi(\alpha-1)} + \mathrm{e}^{\mathrm{i}\pi(\alpha-1)})
 \int_{\epsilon}^{1-\epsilon}x^{\alpha-n-1}f(x)\mathrm{d}x
 \\ 
 = \: & 
 \int_{\epsilon}^{(1+\epsilon)^{-1}}x^{\alpha-n-1}f(x)\mathrm{d}x,
\end{align*}
where we remark that $F(\alpha,1;2-\alpha;z)$ is a single-valued analytic function on the interval $[0,1)$. 
As to the integral on $C_{\epsilon}^{(0)}$, we have
\begin{align*}
 & 
 \frac{1}{2\pi\mathrm{i}}\int_{C_{\epsilon}^{(0)}}z^{-n}f(z)\Psi_{\alpha}(z)\mathrm{d}z
 \\ 
 = \: & 
 - \frac{1}{2\mathrm{i}\sin\pi\alpha}\int_{C_{\epsilon}^{(0)}}
 z^{-n}f(z)(-z)^{\alpha}\mathrm{d}z 
 - 
 \frac{1}{2\pi\mathrm{i}(\alpha-1)}\int_{C_{\epsilon}^{(0)}}
 z^{-n}f(z)F(1-\alpha,1;2-\alpha;z)\mathrm{d}z.
\end{align*}
The first term on the right-hand side is written as
\begin{align*}
 \mbox{(the first term)}
 = \: & 
 \frac{(-1)^{n+1}}{2\mathrm{i}\sin\pi\alpha}
 \int_{C_{\epsilon}^{(0)}}f(z)(-z)^{\alpha-n-1}\mathrm{d}z
 \\ 
 = \: & 
 \frac{(-1)^{n+1}}{2\mathrm{i}\sin\pi\alpha}\int_0^{2\pi}
 f(\epsilon\mathrm{e}^{\mathrm{i}\theta})\{\epsilon\mathrm{e}^{\mathrm{i}(\theta-\pi)}\}^{\alpha-n-1}
 \mathrm{i}\epsilon\mathrm{e}^{\mathrm{i}\theta}\mathrm{d}\theta
 \\ 
 = \: & 
 \frac{\epsilon^{\alpha-n}\mathrm{e}^{-\mathrm{i}\pi\alpha}}{2\sin\pi\alpha}\int_0^{2\pi}
 f(\epsilon\mathrm{e}^{\mathrm{i}\theta})\mathrm{e}^{\mathrm{i}(\alpha-n)\theta}\mathrm{d}\theta
 \\ 
 = \: & 
 \frac{\epsilon^{\alpha-n}\mathrm{e}^{-\mathrm{i}\pi\alpha}}{2\sin\pi\alpha}\int_0^{2\pi}
 \left\{
 \sum_{k=0}^{n-1}\frac{\epsilon^k}{k!}\mathrm{e}^{\mathrm{i}k\theta}f^{(k)}(0)+\mathrm{O}(\epsilon^n)
 \right\}
 \mathrm{e}^{\mathrm{i}(\alpha-n)}
 \mathrm{d}\theta
 \\ 
 = \: & 
 - \sum_{k=0}^{n-1}\frac{\epsilon^{\alpha-n+k}}{k!(n-k-\alpha)}f^{(k)}(0)
 + \mathrm{O}(\epsilon^{\alpha}), 
\end{align*}
and the second term is written as
\begin{align*}
 & \mbox{(the second term)} 
 \\ 
 = \: & 
 - \frac{1}{2\pi\mathrm{i}(\alpha-1)}
 \int_{C_{\epsilon}^{(0)}}z^{-n}\left\{\sum_{k=0}^{\infty}\frac{f^{(k)}(0)}{k!}z^k\right\}
 \left\{\sum_{l=0}^{\infty}\frac{1-\alpha}{1-\alpha+l}z^l\right\}\mathrm{d}z
 \\
 = \: & 
 - \sum_{k+l=n-1}\frac{f^{(k)}(0)}{k!(\alpha-l-1)}
 =  
 - \sum_{k=0}^{n-1}\frac{f^{(k)}(0)}{k!(\alpha-n+k)}.
\end{align*}
As to the integral on $C_{\epsilon}^{(1)}$, from the formula 15.3.10 in \cite{AbramowitzStegun1965}, we have
\begin{gather*}
 \Psi_{\alpha}(z) = 
 z^{-1}\sum_{k=0}^{\infty}\frac{(\alpha)_k}{k!}
 \left\{\psi(k+1) - \psi(k+\alpha) - \log(1 - z^{-1}) \right\}(1 - z^{-1})^k, 
 \intertext{where $\psi(z)$ is the Digamma function: $\psi(z)=\Gamma^{\prime}(z)/\Gamma(z)$ and}
 (\alpha)_0 = 1, \quad 
 (\alpha)_k = \alpha(\alpha+1)(\alpha+2)\cdots(\alpha+k-1) \quad ( \: k = 1, 2, \ldots \: ),
\end{gather*}
and, then, the integral on $C_{\epsilon}^{(1)}$ is of $\mathrm{O}(\epsilon\log\epsilon)$. 
Summarizing the above calculations, we have
\begin{align*}
 \frac{1}{2\pi\mathrm{i}}\oint_C z^{-n}f(z)\Psi_{\alpha}(z)\mathrm{d}z
 = \: & 
 \int_{\epsilon}^{1}x^{\alpha-n-1}f(x)\mathrm{d}x 
 - 
 \sum_{k=0}^{n-1}\frac{\epsilon^{\alpha-n+k}}{k!(n-k-\alpha)}f^{(k)}(0)
 \\
 & 
 - 
 \sum_{k=0}^{n-1}\frac{f^{(k)}(0)}{k!(\alpha-n+k)} 
 + 
 \mathrm{O}(\epsilon\log\epsilon),
\end{align*}
and, taking the limit $\epsilon\downarrow 0$, we have (\ref{eq:complex-integral}).

\hfill\rule{1.5ex}{1.5ex}

\medskip

We can obtain the desired f.p. integral (\ref{eq:fp-integral}) by evaluating 
the complex integral in (\ref{eq:complex-integral}) on the closed integral path 
$C$, which is parameterized by $z = \varphi(u)$, $0\leqq u\leqq u_{\rm p}$, 
by the trapezoidal formula with equal mesh as follows. 
\begin{multline}
 \label{eq:approx-fp-integral}
 \fp\int_0^1 x^{\alpha-1-n}f(x)\mathrm{d}x \simeq I_N^{\alpha,n}[f]
 \\ 
 \equiv
 \frac{h}{2\pi\mathrm{i}}\sum_{k=0}^{N-1}
 \varphi(kh)^{-n}f(\varphi(kh))\Psi_{\alpha}(\varphi(kh))
 \varphi^{\prime}(kh)
 + 
 \sum_{k=0}^{n-1}
 \frac{f^{(k)}(0)}{k!(\alpha-n+k)}
 \\ 
 \left( \: h = \frac{u_{\rm p}}{N} \: \right).
\end{multline}
The hypergeometric function in the definition of $\Psi_{\alpha}(z)$ 
in (\ref{eq:Psi}) is easily evaluated using the continued fraction expansion 
(see \S 12.5 in \cite{Henrici1977}). 
If $C$ is an analytic curve, the complex loop integral is an integral of 
an analytic periodic function on an interval of one period, to which 
the trapezoidal formula with equal mesh is very effective, and, then, 
the approximation formula (\ref{eq:approx-fp-integral}) is very accurate. 
In fact, applying the theorem in \S 4.6.5 in \cite{DavisRabinowitz1984} to 
the approximation of the complex integral in (\ref{eq:approx-fp-integral}), 
we have the following theorem on error estimate of the proposed approximation. 
\begin{thm}
 \label{thm:error-estimate}
 We suppose that
 \begin{itemize}
  \item the parameterization function $\varphi(w)$ of $C$ is analytic 
	in the strip domain
	\begin{equation*}
	 D_d = 
	  \left\{ \: z\in\mathbb{C} \: | \: |\im z| < d \right\} 
	  \quad ( \: d > 0 \: ), 
	\end{equation*}
  \item the domain 
	\begin{equation*}
	 \varphi(D_d) = 
	  \{ \: \varphi(w) \: | \: w\in D_d \: \}
	\end{equation*}
	is contained in $\mathbb{C}\setminus[0,1]$, and
  \item the function $f(z)$ is analytic in $\varphi_{\alpha}(D_d)$.
 \end{itemize}
 Then, we have for arbitrary $0<d^{\prime}<d$ 
 \begin{gather}
  \left| 
   \fp\int_0^1 x^{\alpha-1-n}f(x)\mathrm{d}x - I_N^{(n,\alpha)}[f]
  \right|
  \leqq 
  \frac{u_{\rm p}}{\pi}\mathscr{N}(f,\alpha,n,d^{\prime})
  \frac{\exp(-2\pi d^{\prime} N/u_{\rm p})}{1-\exp(-2\pi d^{\prime}N/u_{\rm p})},
  \intertext{where}
  \mathscr{N}(f,\alpha,n,d^{\prime}) = 
  \max_{|\im w|=d^{\prime}}
  \left|
  \varphi(w)^{-n}f(\varphi(w))\Psi_{\alpha}(w)\varphi^{\prime}(w)
  \right|.
 \end{gather}
\end{thm}
This theorem says that the proposed approximation (\ref{eq:approx-fp-integral}) 
converges exponentially as the number of sampling points $N$ increases 
if $f(x)$ is an analytic periodic function and $C$ is an analytic curve.

We remark here that, if $f(z)$ is real valued on the real axis, we can reduce the number of sampling points by half. 
In fact, in this case, we have
\begin{equation*}
 f(\overline{z}) = \overline{f(z)}
\end{equation*} 
from the reflection principle, and, taking the integral path $C$ to be symmetric with respect to the real axis, 
that is,  
\begin{equation*}
 \varphi(-u) = \overline{\varphi(u)}, \quad \varphi^{\prime}(-u) = - \overline{\varphi^{\prime}(u)}, 
\end{equation*}
we have
\begin{align}
 \nonumber
 & 
 \fp\int_0^1 x^{\alpha-1-n}f(x)\mathrm{d}x \simeq {I^{\prime}}_N^{(n,\alpha)}[f]
 \\ 
 \nonumber
 \equiv \: & 
 \frac{h}{2\pi}\im
 \bigg\{
 \varphi(0)^{-n}f(\varphi(0))\Psi_{\alpha}(\varphi(0))\varphi^{\prime}(0)
 \\ 
 \nonumber 
 & \hspace{15mm}
 + 
 \varphi\left(\frac{u_{\rm p}}{2}\right)^{-n}f\left(\frac{u_{\rm p}}{2}\right)
 \Psi_{\alpha}\left(\varphi\left(\frac{u_{\rm p}}{2}\right)\right)\varphi^{\prime}\left(\frac{u_{\rm p}}{2}\right)
 \bigg\}
 \\ 
 \nonumber
 & 
 + 
 \frac{h}{\pi}\im
 \left\{
 \sum_{k=1}^{N-1}\varphi(kh)^{-n}f(\varphi(kh))\Psi_{\alpha}(\varphi(kh))\varphi^{\prime}(kh)
 \right\}
 + 
 \sum_{k=0}^{n-1}\frac{f^{(k)}(0)}{k!(n-\alpha+k)}
 \\
 \label{eq:approx-fp-integral2}
 & \hspace{80mm}
 \quad \left( \: h = \frac{u_{\rm p}}{2N} \: \right). 
\end{align}
\section{Numerical examples}
\label{sec:example}
\begin{figure}[htbp]
 \begin{center}
  \begin{tabular}{cc}
   \includegraphics[width=0.45\textwidth]{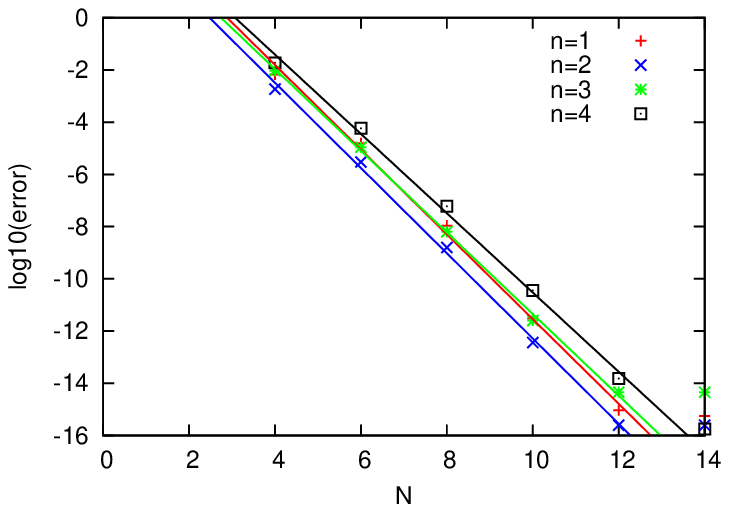} & 
	   \includegraphics[width=0.45\textwidth]{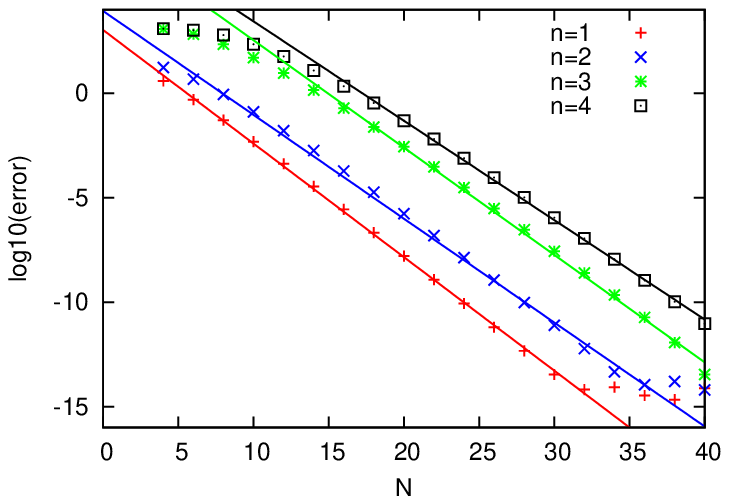}\\
   (i) & (ii)
  \end{tabular}
 \end{center}
 \caption{The relative errors of the proposed method applied 
 to the f.p. integrals (i) and (ii) in (\ref{eq:example}).}
 \label{fig:example}
\end{figure}

In this section, we show some numerical examples which show the effectiveness 
of the proposed method. 
We computed the integrals
\begin{equation}
 \label{eq:example}
 \begin{aligned}
  \mathrm{(i)} \quad & 
  \fp\int_0^1 x^{\alpha-n-1}\mathrm{e}^x\mathrm{d}x = 
  \frac{F(\alpha-n;\alpha+1-n;1)}{\alpha-n}, 
  \\ 
  \mathrm{(ii)} \quad & 
  \fp\int_0^1 \frac{x^{\alpha-n-1}}{1+x^2}\mathrm{d}x  = 
  \frac{1}{\alpha-n}\re F(\alpha-n,1;\alpha+1-n;\mathrm{i})
 \end{aligned}
\end{equation}
with $\alpha=0.1$ 
by the approximation formula (\ref{eq:approx-fp-integral2}). 
All the computations were performed using programs coded in C++ 
with double precision working. 
The complex integral path $C$ was taken as the ellipse 
\begin{equation*}
 C \: : \: 
  z = \frac{1}{2} + \frac{1}{4}\left(\rho+\frac{1}{\rho}\right)\cos u + 
  \frac{1}{4}\left(\rho-\frac{1}{\rho}\right)\sin u, 
  \quad 0 \leqq u \leqq u_{\rm p} \quad ( \: \rho > 0 \: )
\end{equation*}
with $\rho=10$ for the integral (i) and $\rho=2$ for the integral (ii). 
Figure \ref{fig:example} show the relative errors of the proposed method 
applied to the integrals (i) and (ii) as functions of the number of sampling 
points $N$. 
From these figures, the errors of the proposed method decays exponentially 
as $N$ increases, and the decay rate of the error does not depend much on $n$. 
Table \ref{tab:error} shows the decay rates of the errors of the proposed method 
applied to the f.p. integrals (i) and (ii). 
\begin{table}[htbp]
 \caption{The decay rates of the errors of the proposed method allied to 
 the the f.p. integrals (i) and (ii) in (\ref{eq:example}).}
 \begin{center}
  \begin{tabular}{ccccc}
   \hline
   $n$ & 1 & 2 & 3 & 4 \\
   \hline
   integral (1) & $\mathrm{O}(0.024^N)$ & $\mathrm{O}(0.023^N)$ & 
			   $\mathrm{O}(0.027^N)$ & $\mathrm{O}(0.030^N)$ 
				   \\ 
   \hline
   integral (2) & $\mathrm{O}(0.28^N)$ & $\mathrm{O}(0.32^N)$ & 
			   $\mathrm{O}(0.31^N)$ & $\mathrm{O}(0.33^N)$
	   \\
   \hline
  \end{tabular}
 \end{center}
 \label{tab:error}
\end{table}
\section{Summary}
\label{sec:summary}
In this paper, we proposed a numerical method of computing Hadamard finite part 
integrals with a non-integral power singularity on an endpoint. 
In the proposed method, we express the desired f.p. integral using a complex 
loop integral, and obtain the f.p. integral by evaluating the complex integral 
by the trapezoidal formula with equal mesh. 
Theoretical error estimate and some numerical examples showed the exponential 
convergence of the proposed method. 

We can obtain similarly f.p. integrals on an infinite interval. 
This will be reported in a future paper. 
\bibliographystyle{plain}
\bibliography{arxiv2019_3}

\begin{thebibliography}{10}

\bibitem{AbramowitzStegun1965}
M.~Abramowitz and Irene A.~Stegun (eds.).
\newblock {\em Handbook of Mathematical Functions with Formulas, Graphs and
  Mathematical Tables}.
\newblock Dover, New York, 1965.

\bibitem{Bialecki1990a}
B.~Bialecki.
\newblock A sinc-hunter quadrature rule for cauchy principal value integrals.
\newblock {\em Math. Comput.}, 55:665--681, 1990.

\bibitem{Bialecki1990b}
B.~Bialecki.
\newblock A sinc quadrature rule for hadamard finite-part integrals.
\newblock {\em Numer. Math.}, 57:263--269, 1990.

\bibitem{DavisRabinowitz1984}
P.~J. Davis and P.~Rabinowitz.
\newblock {\em Methods of Numerical Integration, Second Ed.}
\newblock Academic Press, San Diego, 1984.

\bibitem{ElliotPaget1979}
D.~Elliot and D.~F. Paget.
\newblock Gauss type quadrature rules for cauchy principal value integrals.
\newblock {\em Math. Comput.}, 33:301--309, 1979.

\bibitem{EstradaKanwal1989}
R.~Estrada and R.~P. Kanwal.
\newblock Regularization, pseudofunction, and hadamard finite part.
\newblock {\em J. Math. Anal. Appl.}, 141:195--207, 1989.

\bibitem{Henrici1977}
P.~Henrici.
\newblock {\em Applied and Computational Complex Analysis}, volume~2.
\newblock John Wiley \& Sons, New York, 1977.

\bibitem{Ogata2019b}
H.~Ogata.
\newblock A numerical method for hadamard finite-part integrals with an
  integral power singularity at an endpoint, 2019.
\newblock arXiv:1909.08872v1 [math.NA].

\bibitem{OgataHirayama2018}
H.~Ogata and H.~Hirayama.
\newblock Numerical integration based on hyperfunction theory.
\newblock {\em J. Comput. Appl. Math.}, 327:243--259, 2018.

\bibitem{OgataSugiharaMori2000}
H.~Ogata, M.~Sugihara, and M.~Mori.
\newblock De-type quadrature formulae for cauchy principal-value integrals and
  for hadamard finite-part itnegrals.
\newblock In {\em Proceedings of the Second ISAAC Congress}, volume~1, pages
  357--366, 2000.

\bibitem{Paget1981}
D.~F. Paget.
\newblock The numerical evaluation of hadamard finite-part integrals.
\newblock {\em Numer. Math.}, 36:447--453, 1981.

\bibitem{TakahasiMori1974}
H.~Takahasi and M.~Mori.
\newblock Double exponential formulas for numerical integration.
\newblock {\em Publ. Res. Inst. Math. Sci., Kyoto Univ.}, 339:721--741, 1978.

\end{thebibliography}
\end{document}